\documentclass{elsarticle}
\usepackage{amssymb, color}
\usepackage{amsmath}
 \usepackage{pifont}
 \usepackage{caption}
 \usepackage[utf8]{inputenc}
\usepackage{setspace}
\makeatletter
\def\LaTeX{\leavevmode L\raise.42ex
    \hbox{\kern-.3em\size{\sf@size}{0pt}\selectfont A}\kern-.15em\TeX}
\makeatother

\newcommand{\BibTeX}{{\rm B\kern-.05em{\sc
          i\kern-.025emb}\kern-.08em\TeX}}

\makeatletter
\def\@currentlabel{2.1}\label{e:dispaa}
\def\@currentlabel{2.21}\label{e:dispau}
\def\@currentlabel{2.22}\label{e:dispav}
\def\@currentlabel{2.23}\label{e:dispaw}
\def\@currentlabel{2.24}\label{e:dispax}
\def\theequation{\thesection.\@arabic\c@equation}
\makeatother

\renewcommand{\theequation}{\arabic{section}.\arabic{equation}}

\newcommand{\R}{\mathbb R}

\newcommand{\e}{\epsilon}

\def \D{\Delta}
\def \O{\Omega}
\newtheorem{thm}{Theorem} [section]
\newtheorem{lem}{Lemma} [section]
\newtheorem{prop}{Proposition} [section]

\newtheorem{taggedtheoremx}{Theorem}
\newenvironment{taggedtheorem}[1]
 {\renewcommand\thetaggedtheoremx{#1}\taggedtheoremx}
 {\endtaggedtheoremx}

\renewcommand{\theequation}{\thesection.\arabic{equation}}
\renewcommand{\thesection}{\arabic{section}}
\renewcommand{\theequation}{\thesection.\arabic{equation}}
\let\ssection=\section\renewcommand{\section}{\setcounter{equation}{0}\ssection}

\def \p{\partial}
\def \wt{\widetilde}

\begin{document}

\begin{frontmatter}

\title{Liouville theorems for stable at infinity solutions of Lane-Emden system}
\author[fd]{Foued Mtiri}
\ead{mtirifoued@yahoo.fr}
\author[ah]{Dong Ye\corref{cor1}}
\ead{dong.ye@univ-lorraine.fr}
\address[fd]{ANLIG, UR13ES32, University of Tunis El-Manar, 2092 El Manar II, Tunisia}
\address[ah]{IECL, UMR 7502, Universit\'{e} de Lorraine, 3 rue Augustin Fresnel, 57073 Metz, France}
\begin{abstract}
We consider the Lane-Emden system $-\Delta u = v^p$, $-\Delta v= u^\theta$ in $\mathbb{R}^N$, and we prove the nonexistence of smooth positive solutions which are stable outside a compact set,
for any $p, \theta > 0$ under the Sobolev hyperbola.
\end{abstract}
\end{frontmatter}
 \section{Introduction}
\setcounter{equation}{0}

Consider the classical Lane-Emden system
\begin{align}\label{1.1}
-\Delta u = v^p, \quad-\Delta v= u^\theta,\quad u,v>0\quad\mbox{in }\; \mathbb{R}^N, \quad\mbox{where }\; p,\theta >0.
\end{align}
There is a famous conjecture who states that: {\sl Let $p, \theta > 0$. If the pair $(p, \theta)$ is subcritical, i.e.~if
\begin{align}\label{LE}
\frac{1}{p+1} + \frac{1}{\theta +1} > \frac{N-2}{N},
\end{align}
then there is no smooth solution to \eqref{1.1}.}

\medskip
The critical curve given by the equality in \eqref{LE} is called the Sobolev hyperbola, which is introduced independently by Mitidieri \cite{em} and
Van der Vorst \cite{van}, it plays a crucial role in the analysis of \eqref{1.1}. It is well known that if $(p, \theta)$ lies on or above the Sobolev hyperbola,
\eqref{1.1} admits radial classical solutions (see \cite{em1, Zs}), and the Lane–Emden conjecture can be restated as the following: There has
no smooth solution to \eqref{1.1} if the positive pair $(p, \theta)$ lies below the Sobolev critical hyperbola.

\medskip
The conjecture is proved to be true for radial functions by Mitidieri \cite{em1}, Serrin-Zou \cite{Zs1}. For the full conjecture, Souto \cite{s}, Mitidieri \cite{em1} and Serrin-Zou \cite{Zs} proved that there is no supersolution to \eqref{1.1}, if $p\theta \leq 1$ or $\max(\alpha, \beta) \geq N-2$, where
 \begin{align}\label{ab}
  \alpha = \frac{2(p+1)}{p\theta - 1}, \quad \beta = \frac{2(\theta +1)}{p\theta - 1}, \quad p\theta > 1.
 \end{align}
 Moreover,
we can check readily that if $p\theta > 1$, the condition \eqref{LE} is equivalent to
\begin{align}
 \label{LEbis} N < 2 + \alpha + \beta.
\end{align}
Therefore, the Lane-Emden conjecture is true in dimensions $N = 1, 2$. More recently, the conjecture is proved in dimensions $N = 3, 4$, by Souplet and his collaborators, see \cite{pqs, sou}.
For $N \geq 5$, the conjecture is known to be true for $(p, \theta)$ verifying \eqref{LE} and one of the following extra conditions:
\begin{itemize}
 \item If $p, \theta < \frac{N+2}{N-2}$, see Felmer-de Figuereido \cite{ff}.
 \item If $\max(p, \theta) \geq N-3$, see Souplet \cite{sou}.
 \item If $\min(\alpha, \beta) \geq \frac{N-2}{2}$, see Busca-Man\'asevich \cite{bm}\footnote{In \cite{bm}, there is another extra condition, which is no longer necessary after the work of Souplet \cite{sou}.}.
 \item If $p = 1$ or $\theta = 1$, see Lin \cite{lin}.
\end{itemize}
These partial results enable us a more restrictive new region for the exponents $(p, \theta)$, which is illustrated by the following figure. In other words, the Lane-Emden conjecture stands open for $N \geq 5$, $p,\theta > 0$ such that
\begin{align*}
p, \theta \ne 1, \quad \min(\alpha, \beta) < \frac{N-2}{2} \quad\mbox{and}\quad \max(\alpha, \beta) < N-3.
\end{align*}
\begin{figure}[h]
\begin{center}
\includegraphics[width=8cm,height=6cm]{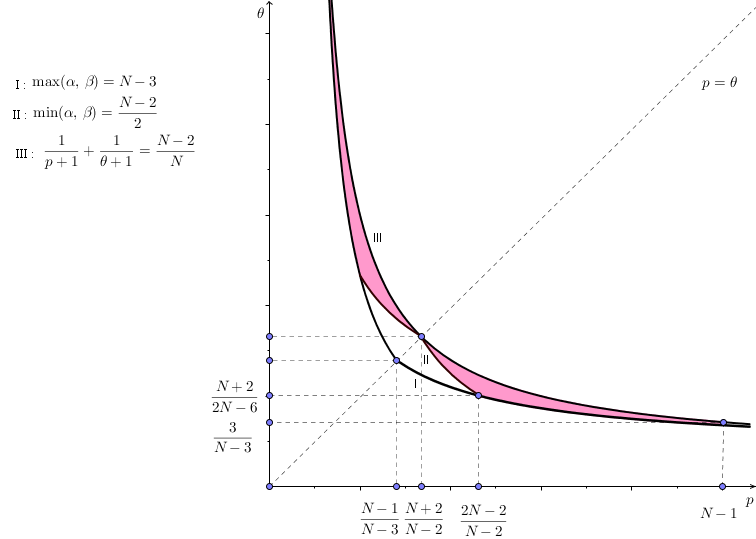}
\caption[le titre]{The remained open region (shaded) for
the Lane-Emden conjecture ($N \geq 5$)}
\label{monlabel}
\end{center}
\end{figure}

On the other hand, in the last decade, many efforts were made to obtain some Liouville type result for solutions with finite Morse index, or more generally,
which are stable at infinity. To define the
notion of stability, we consider a general system given by
\begin{align}\label{1.222}
 -\Delta u = f(x,v),\quad -\Delta v= g(x,u)\;\; \mbox{in $K$, a bounded regular domain }\; \subset \R^N,
\end{align}
where $f,g\in C^1(K \times \mathbb{R}).$  Following Montenegro \cite{MO}, a smooth solution $(u,v)$ of \eqref{1.222} is said to
be stable in $K$ if the following eigenvalue problem
\begin{align*}
 -\Delta \xi = f_{v}(x,v)\zeta +\eta \xi, \quad -\Delta \zeta = g_{u}(x,u)\xi + \eta\zeta \quad \mbox{in }\, K
\end{align*}
has a nonnegative eigenvalue $\eta$, with a positive smooth eigenfunctions pair $(\xi, \zeta)$. We say that a pair of solutions $(u, v)$ to \eqref{1.1} is stable outside a compact
set or stable at infinity, if there is a compact set $K\subset \R^N$ such that $(u, v)$ is stable in any bounded domain of $\R^N\backslash K$.

\medskip
For the corresponding second order equation
\begin{align}
\label{2nd}
\Delta u+ |u|^{q-1}u =0 \quad \mbox{in } \R^N, \quad q > 1,
\end{align}
Farina \cite{Far} obtained the optimal Liouville type result for solutions stable at infinity. Indeed, he proved that a smooth nontrivial solution to \eqref{2nd} exists,
if and only if $q> p_{JL}$ and $N\geq 11,$ or $q =\frac{N+2}{N-2}$ and $N\geq 3.$ Here $p_{JL}$ denotes the so-called Joseph-Lundgren exponent (see \cite{CW, Far}).
For the biharmornic equation $\Delta^2 u = |u|^{q-1}u$, $q > 1$, D\'{a}vila-Dupaigne-Wang-Wei \cite{ddww} derived a striking monotonicity formula,
which led them to the optimal classification result for solutions stable at infinity, using blow down analysis.

\medskip
Coming back to the Lane-Emden system \eqref{1.1}, Chen-Dupaigne-Ghergu \cite{WLD} studied the stability of radial solutions when $p, \theta \geq 1.$
They introduced a new critical hyperbola, called the Joseph-Lundgren curve.  More precisely, they proved that if $p, \theta \geq 1,$  then a radial solution
of \eqref{1.1} is unstable if and only if $N \leq 10$, or $N \geq 11$ and
\begin{equation*}
  \left[\frac{(N-2)^2-(\alpha-\beta)^2}{4}\right]^2<p\theta\alpha\beta(N-2-\alpha)(N-2-\beta).
\end{equation*}

Moreover, Cowan proved in \cite{cow} that if $p, \theta\geq 2$ and $N \leq 10$, there does not exist any stable solution (radial or not) to \eqref{1.1}.
Recently, Hajlaoui-Harrabi-Mtiri \cite{Hfh} established some Liouville theorems for smooth stable solutions of \eqref{1.1} with $p > 1$, see Theorem {\bf A} below.
We mention also the celebrated result of Ramos \cite{rm}, which states that if $p, \theta > 1$ satisfies \eqref{LE}, then the system
$$-\Delta u = |v|^{p-1}v, \;\; -\Delta v = |u|^{\theta-1}u \quad \mbox{in } \R^N$$ does not admit any smooth solutions
having finite {\sl relative} Morse index in the sense of Abbondandolo.

\medskip
In this paper, our motivation are twofold. We want to obtain the classification results for solutions (radial or not) to \eqref{1.1} which are just stable at infinity,
and we want to handle the case where $p, \theta$ are allowed to be less than 1.
So a natural question is: Can we prove the Lane-Emden conjecture with the extra condition that $(u, v)$ is
stable at infinity? The answer is affirmative.
\begin{thm}
\label{main4}
For any $p, \theta > 0$ satisfying \eqref{LE}, the system \eqref{1.1} has no classical solution stable outside a compact set.
\end{thm}

If $\theta = p$, using Souplet's comparison result (Lemma 2.7 in \cite{sou}), we get $u \equiv v$, so the optimal classification result for solutions stable
at infinity was already given by Farina. The classification is also known for $p\theta \leq 1$ as mentioned above. Without loss of generality, we consider only $\theta > p > 0$ and $p \theta > 1$. As we will see soon, the $\theta > p \geq 1$ case can be handled by
the results in \cite{Hfh}, so our main concern is the case
$$\theta p > 1 > p > 0.$$

\medskip
Let $(u, v)$ be a smooth solution to \eqref{1.1} with $\theta > p^{-1} > 1 > p > 0$. Our approach is based on the formal equivalence noticed in \cite{DEN, dos}, between
the Lane-Emden system \eqref{1.1} and a fourth order problem, called the $m$-biharmonic equation. More precisely, let $m:=\frac{1}{p}+1>2$, as $v =(-\Delta u)^{m-1}$, we derive
that $u$ satisfies $\Delta^{2}_{m} u := \Delta (|\Delta u|^{m-2}\Delta u) = u^\theta$ in $\R^N$. So we are led to consider $\theta > m-1 > 1$ and
\begin{equation}\label{1}
\Delta^{2}_{m} u := \Delta (|\Delta u|^{m-2}\Delta u) =|u|^{\theta-1}u.
\end{equation}

\smallskip
Let $\Omega \subset \R^N$, we say that $u \in W^{2,m}_{loc}(\Omega)\cap L^{\theta +1}_{loc}(\Omega)$ is a weak solutions of \eqref{1} in $\Omega$, if for any regular
bounded domain $K \subset \Omega$,
$u$ is a critical point of the following functional
$$I(v)=\frac{1}{m}\int_K |\D v|^m dx-\frac{1}{\theta+1}\int_K |v|^{\theta+1} dx, \quad \forall\; v \in W^{2,m}(K)\cap L^{\theta +1}(K).$$
Naturally, a weak solution to \eqref{1} is said stable in $\Omega \subset \R^N$, if
\begin{equation}\label{fb1}
\Lambda_u (h):= (m-1)\int_{\Omega} |\D u|^{m-2}|\D h|^2 dx-\theta \int_{\Omega}|u|^{\theta-1} h^2dx\geq0,\;\;\forall\; h\in C_c^2(\Omega).
\end{equation}
A key point for our approach is to remark a relationship between the stability for the system \eqref{1.1} and the stability for the equation \eqref{1} (see Lemma 2.1 below). This will permit us to handle the case $0 < p  < 1$ in \eqref{1.1} by using the structure of the $m$-biharmonic equation. In fact, we can prove the following Liouville type result.
\begin{thm}\label{main3}
Let $\theta > m-1 > 1$ and $u\in W^{2,m}_{loc}(\mathbb{R}^N)\cap L^{\theta + 1}_{loc}(\R^N)$ be a weak solution of \eqref{1} which is stable outside a compact set. Assume that
\begin{align}\label{new8}
N < \frac{2m(\theta + 1)}{\theta -(m-1)},
 \end{align}
then $u\equiv 0.$
\end{thm}
A direct calculation yields that if $p\theta > 1$ (or equivalently $\theta > m-1$),
\begin{align*}
 N < 2 + \alpha + \beta = \frac{2(p+1)(\theta+1)}{p\theta - 1} \;\; \Leftrightarrow \;\; \eqref{new8} \;\; \Leftrightarrow \;\; \theta < \frac{N(m-1)+2m}{(N-2m)_+}.
\end{align*}
It means that the range of pairs $(p, \theta)$ satisfying \eqref{LE} and $p\theta > 1$ corresponds exactly to the subcritical case of the $m$-biharmonic equation \eqref{1}.

\medskip
Another crucial step in our approach is to classify first the stable solutions of \eqref{1.1}, see also Proposition \ref{p12bis} below for the $m$-biharmonic equation.
\begin{prop}\label{p12}
If $p, \theta > 0$ satisfies \eqref{LE}, then \eqref{1.1} has no smooth stable solution.
\end{prop}

Establishing a Liouville type result for stable solution of \eqref{1.1} or \eqref{1} is delicate, even we can borrow some ideas from \cite{wy, ddww}. We use as usual the stability to
get some integral estimates, but the integrations by parts argument yields here many terms which are difficult to control,
for example, the local $L^m$ norm of $\nabla u$, see Lemma \ref{l.2.7a} below. Furthermore, the classification of weak solutions stable at infinity to \eqref{1} is more involved
than to handle \eqref{1.1}, since the weak solutions to \eqref{1} are not $C^2$ functions. We will derive a variant of the Pohozaev identity with cut-off functions,
which allows us to avoid the spherical integral terms in the standard Pohozaev identity.

\medskip
The paper is organized as follows. In section 2, we give the proof
of Proposition \ref{p12}. The proofs  of Theorem \ref{main4} and Theorem \ref{main3} are given  respectively in sections 3 and 4. In the following, $C$
denotes always a generic positive constant, which could be changed from one line to another.

\section{Classification of stable solutions}
\setcounter{equation}{0}
We prove here Proposition \ref{p12}. As mentioned before, we need only to consider the case $\theta > p$ and $p \theta > 1$. We split the proof into two cases:
$\theta > p \geq 1$ and $\theta > p^{-1} > 1 > p > 0$.

\subsection{ The case $\theta>p \geq 1$.}
Let us recall a consequence of Theorem 1.1 (with $\alpha = 0$ there) in \cite{Hfh}.
\begin{taggedtheorem}{A} Let $x_0$ be the largest root of the polynomial
\begin{align}\label{newH}
H(x)=x^4 -p\theta\alpha\beta \left[4x^{2}-2(\alpha+\beta)x+1\right].
\end{align}
  \begin{enumerate}
\item[(i)]
 If $\frac{4}{3}< p \leq \theta$ then \eqref{1.1} has no stable classical solution if $N<2+2x_0.$
\item[(ii)] If $1 \leq p\leq \min(\frac{4}{3}, \theta)$ and $p\theta > 1$, then \eqref{1.1} has no stable classical solution, if
$$N < 2 + 2x_0 \left[\frac{p}{2}+\frac{(2-p)(p \theta -1)}{(\theta+p-2)(\theta+1)}\right].$$
\end{enumerate}
\end{taggedtheorem}
Performing the change of variables $x=\frac{\beta}{2}s$ in \eqref{newH}, a direct computation shows that $H(x)=\left(\frac{\beta}{2}\right)^4L(s)$ where
$$L(s):=s^4-\frac{16p\theta(p+1)}{\theta+1}s^2+\frac{16p\theta(p+1)(p+\theta+2)}{(\theta+1)^2}s-\frac{16p\theta(p+1)^2}{(\theta+1)^2}.$$
Denote by $s_0$  the largest root of $L,$ hence $x_0=\frac{\beta}{2}s_0$ and $ H(x)<0$ if and only if $L(s)<0$. For $\theta > p \geq 1$, there holds
\begin{align*}
L(p+1) & = (p+1)^4 -\frac{16p\theta(p+1)^3}{(\theta+1)} +\frac{16p\theta(p+1)^2(p+\theta+2)}{(\theta+1)^2} -\frac{16p\theta(p+1)^2}{(\theta+1)^2}\\
& = (p+1)^4 -\frac{16p\theta(p+1)^3}{(\theta+1)} + \frac{16p\theta(p+1)^2}{(\theta+1)} + \frac{16p\theta(p+1)^3}{(\theta+1)^2}
-\frac{16p\theta(p+1)^2}{(\theta+1)^2}\\
& = (p+1)^2\left[(p+1)^2 -\frac{16p^{2}\theta}{(\theta+1)} + \frac{16p^{2}\theta}{(\theta+1)^2}\right]\\
&=\left(\frac{p+1}{\theta+1}\right)^{2}\left[(p+1)^2(\theta+1)^2 -16p^{2}\theta^{2}\right] < 0.
\end{align*}
The last inequality holds true since
$$4p\theta- (p+1)(\theta + 1) > 4p^2 - (p+1)^2 \geq 0, \quad \forall\; \theta > p \geq 1.$$
As $\lim_{s\rightarrow \infty}L(s)= \infty,$  it follows that $s_0>p+1.$ We get then
 \begin{align*}
2x_0>(p+1)\beta=2+\alpha+\beta,\quad \forall\; \theta > p\geq 1.
\end{align*}
If $p>\frac{4}{3}$, by $(i)$ of Theorem \textbf{A}, the system \eqref{1.1} has no classical stable solution if  $N<2+\alpha+\beta$. Suppose now $1\leq p \leq \min(\frac{4}{3}, \theta)$.
Observe that for all $\theta \geq p \geq 1$,
     \begin{align*}
\left[p+\frac{2(2-p)(p \theta -1)}{(\theta+p-2)(\theta+1)}\right]\beta \geq \alpha + \beta & \Leftrightarrow \left[p+\frac{2(2-p)(p \theta -1)}{(\theta+p-2)(\theta+1)}\right](\theta + 1) \geq p + \theta + 2\\
& \Leftrightarrow p\theta - 1 + \frac{2(2-p)(p \theta -1)}{\theta+p-2} \geq \theta + 1\\
& \Leftrightarrow (p\theta - 1)\left[ 1 + \frac{2(2-p)}{\theta+p-2}\right] \geq \theta + 1\\
& \Leftrightarrow (p\theta - 1)(\theta + 2 - p) \geq (\theta+p-2)(\theta + 1)\\
& \Leftrightarrow p\theta^2 - \theta + (2 - p)p\theta \geq \theta^2 + (p-1)\theta\\
& \Leftrightarrow (p-1)(\theta - p) \geq 0.
\end{align*}
As $s_0>p+1\geq 2,$ we have $x_0 = \frac{\beta s_0}{2} \geq \beta$ and
$$2+\alpha+\beta \leq 2 + \beta \left[p+\frac{2(2-p)(p \theta -1)}{(\theta+p-2)(\theta+1)}\right] \leq 2 + x_0\left[p+\frac{2(2-p)(p \theta -1)}{(\theta+p-2)(\theta+1)}\right].$$
If $N < 2+\alpha+\beta,$ using $(ii)$ of Theorem \textbf{A}, we are done.

\medskip
To conclude, for all $\theta > p \geq 1$ and $N < 2+\alpha+\beta,$ \eqref{1.1} has no smooth stable solution. \qed

\subsection{The case $\theta p > 1>p > 0$. }
Here we handle the case  $0<p<1.$ First of all, we need the following lemma which plays an important role in dealing with  Proposition \ref{p12}.
\begin{lem}
\label{l.1} Let $(u,v)$ be a solution of system \eqref{1.1} with $\theta > \frac{1}{p}:=m-1>1$. Suppose that $(u, v)$ is stable in a regular bounded domain $\Omega$,
then $u$ is a stable solution of equation \eqref{1}.
\end{lem}
\noindent{\bf Proof.} By the definition of stability, there exist smooth positive functions $\xi$, $\zeta$ and $\eta \geq 0$ such that
$$-\Delta \xi = pv^{p-1}\zeta + \eta \xi, \; \; -\Delta\zeta = \theta u^{\theta -1}\xi + \eta\zeta \quad \mbox{in }\; \Omega.$$
Using $(\xi, \zeta)$ as super-solution, $(\min_{\overline\Omega}\xi, \min_{\overline\Omega}\zeta)$ as sub-solution, and the standard monotone iterations, we can claim that
there exist positive smooth functions $\varphi$, $\chi$ verifying
\begin{align*}
 -\Delta \varphi = p v^{p-1}\chi, \quad -\Delta \chi = \theta u^{\theta-1}\varphi\quad \mbox{in }\, \Omega.
\end{align*}
Therefore, we have
\begin{align*}
\theta u^{\theta-1}\varphi=\Delta\left(\frac{1}{p} v^{1-p} \Delta \varphi\right) \quad \mbox{in}\;\;\Omega.
\end{align*}
Let $\gamma \in C_c^2(\Omega)$. Multiplying the above equation by $\gamma^{2}\varphi^{-1}$ and
integrating by parts, there holds
\begin{align}\label{0.255}
\begin{split}
 \int_{\Omega} \theta u^{\theta-1}\gamma^{2}dx & = \frac{1}{p}\int_{\Omega} v^{1-p} \Delta \varphi\Delta(\gamma^{2}\varphi^{-1})dx\\
&= \frac{1}{p} \int_{\Omega}v^{1-p} \Delta \varphi\left(-4\gamma\frac{\nabla \varphi\cdot\nabla\gamma}{\varphi^{2}}+\frac{2|\nabla\gamma |^2}{\varphi}+\frac{2\gamma\D\gamma}{\varphi}+\frac{2\gamma^{2}|\nabla\varphi |^2}{\varphi^{3}} -\frac{\gamma^{2}\D\varphi}{\varphi^{2}}\right) dx.
\end{split}
\end{align}
Using Cauchy-Schwarz's inequality and the fact that $-\Delta \varphi >0,$ we get
\begin{align}\label{2.La}
\begin{split}
\left| -4\int_{\Omega}\frac{v^{1-p}}{p} \Delta \varphi\frac{\nabla \varphi\cdot\nabla\gamma}{\varphi^{2}}\gamma dx\right|
 \leq -2\int_{\Omega}\frac{v^{1-p}}{p} \Delta \varphi\frac{|\nabla\gamma |^2}{\varphi} dx -2\int_{\Omega}\frac{v^{1-p}}{p} \Delta \varphi\frac{\gamma^{2}|\nabla\varphi |^2}{\varphi^{3}} dx.
\end{split}
\end{align}
Combining \eqref{0.255} and  \eqref{2.La}, one obtains, using again the Cauchy-Schwartz inequality,
\begin{align*}
 \int_{\Omega} \theta u^{\theta-1}\gamma^{2}dx &
\leq \frac{2}{p} \int_{\Omega}v^{1-p} \Delta \varphi\frac{\gamma\D\gamma}{\varphi} dx-\frac{1}{p} \int_{\Omega}v^{1-p}\frac{(\Delta \varphi)^{2}}{\varphi^{2}}\gamma^{2} dx\\
& \leq  \frac{1}{p}\int_{\Omega}v^{1-p}\frac{(\Delta \varphi)^{2}}{\varphi^{2}}\gamma^{2} dx + \frac{1}{p}\int_{\Omega}v^{1-p}(\D\gamma)^{2} dx -\frac{1}{p} \int_{\Omega}v^{1-p}\frac{(\Delta \varphi)^{2}}{\varphi^{2}}\gamma^{2} dx\\
& = \frac{1}{p} \int_{\Omega}v^{1-p}(\D\gamma)^{2} dx.
\end{align*}
Recall that $p = \frac{1}{m-1}$ and $(-\Delta u )^{\frac{1}{p}}= v,$ we obtain the desired result \eqref{fb1}. \qed

\medskip
Therefore, to prove Proposition \ref{p12} in the case $p \in (0, 1)$ and $p\theta > 1$, we need only to prove
\begin{prop}\label{p12bis}
Let $\theta > m - 1 > 1$, if $u$ is a weak stable solution to the equation \eqref{1} in $\R^N$ with $N$ verifying \eqref{new8}, then $u \equiv 0$.
\end{prop}

To prove Proposition \ref{p12bis}, we use first the stability condition \eqref{fb1} to get the  following crucial lemma which provides an important integral estimate for $u$ and $\D u$.
 \begin{lem}
\label{lemnewBN} Let $u\in W^{2,m}_{loc}(\Omega)\cap L^{\theta+1}_{loc}(\Omega)$ be a weak stable solution of \eqref{1} in $\Omega$, with $\theta>m-1 > 1$. Then, for any integer  $$k\geq \max \left( m, \frac{m(\theta+1)}{2(\theta+1-m)}\right) ,$$
there exists a positive constant  $ C=C(N, \epsilon, m, k)$ such that for any $\zeta\in C_c^2(\Omega)$ satisfying $0 \leq\zeta \leq 1$,
\begin{align}\label{new7}
\begin{split}
\int_{\Omega}|\D u|^{m} \zeta^{4k} dx +\int_{\Omega} |u|^{\theta+1}\zeta^{4k} dx
\leq C\left[\int_{\Omega}\left(|\D \zeta|^{m}+|\nabla \zeta|^{2m}
 +|\nabla^{2} \zeta|^{m}\right)^{\frac{\theta+1}{\theta-(m-1)}} dx\right].
 \end{split}
\end{align}
\end{lem}
\noindent{\bf Proof.} For any $\epsilon \in (0, 1)$ and $\eta \in C^2(\Omega)$, there holds
\begin{align}
\label{newesTt47}
  \begin{split}
\int_{\Omega}|\Delta u|^{m-2} [\Delta (u\eta)]^2 dx = & \; \int_{\Omega}|\Delta u|^{m-2}\left(u \D \eta+2\nabla u\nabla \eta + \eta\D u\right)^{2} dx\\
\leq & \;\left(1+\epsilon\right)\int_{\Omega}|\D u|^m {\eta}^2 dx + \frac{C}{\e}\int_{\Omega}|\Delta u|^{m-2}
  \Big(u^2|\D\eta|^2 +  |\nabla u|^2|\nabla \eta|^2\Big)dx.
    \end{split}
   \end{align}
Take $\eta = \zeta^{2k}$ with $\zeta \in C_c^2(\Omega)$, $0 \leq\zeta \leq 1$ and $k \geq m > 2$. Apply Young's inequality, we get
     \begin{align*}
  \int_{\Omega}|u|^{2}|\Delta u|^{m-2}|\D (\zeta^{2k})|^2 dx & \leq C_{k}\int_{\Omega}|u|^{2}|\Delta u|^{m-2}\left(|\D \zeta|^{2}+|\nabla \zeta|^4\right)\zeta^{4k-4} dx\\
& \leq \epsilon^{2} \int_{\Omega}|\D u|^m \zeta^{4k} dx +
  C_{\epsilon, k, m}\int_{\Omega} |u|^m \left(|\D \zeta|^{2}+|\nabla \zeta|^4\right)^{\frac{m}{2}}\zeta^{4k-2m} dx
\end{align*}
and
  \begin{align*}
  \int_{\Omega}|\Delta u|^{m-2} |\nabla u|^2|\nabla (\zeta^{2k})|^2 dx & =4 k^{2}\int_{\Omega}|\Delta u|^{m-2} |\nabla u|^2|\nabla \zeta|^2\zeta^{4k-2} dx\\
& \leq  \epsilon^{2} \int_{\Omega}|\D u|^m \zeta^{4k} dx +
  \frac{C_{m,k}}{\epsilon^{2}}\int_{\Omega} |\nabla u|^m |\nabla \zeta|^{m}\zeta^{4k-m} dx.
\end{align*}
Inserting the two above estimates into \eqref{newesTt47}, we arrive at
\begin{align}\label{newest4}
  \begin{split}
\int_{\Omega}|\Delta u|^{m-2} [\Delta (u\zeta^{2k})]^2 dx
 \leq & \; \left(1+C\epsilon\right) \int_{\Omega}|\D u|^m \zeta^{4k} dx+\frac{C_{m,k}}{\epsilon^{3}}\int_{\Omega} |\nabla u|^m |\nabla \zeta|^{m}\zeta^{4k-m} dx\\
  &+ C_{\epsilon, m, k}\int_{\Omega} |u|^m \left(|\D \zeta|^{2}+|\nabla \zeta|^4\right)^{\frac{m}{2}}\zeta^{4k-2m} dx.
  \end{split}
   \end{align}

\medskip
We need also the following technical lemma, which proof is given later.
 \begin{lem}\label{l.2.7a}
Let $k \geq m/2 > 1$ and $\epsilon > 0$, there exists $C_{N, \e, m, k} >0$ such that for any $u\in W^{2,m}_{loc}(\Omega)$ verifying \eqref{fb1} and $\zeta \in C_c^{\infty}(\Omega)$ with $0 \leq\zeta \leq 1$, there holds
 \begin{align}
 \label{new13}
\int_{\Omega} |\nabla u|^m |\nabla \zeta|^{m}\zeta^{4k-m} dx \leq \epsilon\int_{\Omega}|\D u|^{m}\zeta^{4k} dx +
C_{N, \epsilon, m ,k} \int_{\Omega}|u|^{m}\left(|\nabla \zeta|^{2m} +|\nabla^{2} \zeta|^{m}\right)\zeta^{4k-2m}dx.
\end{align}
\end{lem}

Using Lemma \ref{l.2.7a} with $\epsilon^4$ and \eqref{newest4}, we see that
\begin{align}\label{nt4}
  \begin{split}
\int_{\Omega}|\Delta u|^{m-2} [\Delta (u\zeta^{2k})]^2 dx
  \leq& \; C_{N, \epsilon, m,k}\int_{\Omega} |u|^{m} \left(|\D \zeta|^{m}+|\nabla \zeta|^{2m} +|\nabla^{2} \zeta|^{m}\right)\zeta^{4k-2m} dx\\
  & \;+\left(1 + C_{m,k}\epsilon\right) \int_{\Omega}|\D u|^{m} \zeta^{4k} dx.
  \end{split}
   \end{align}
Thanks to the approximation argument, the stability property \eqref{fb1} holds true with $u\zeta^{2k}$. We deduce then, for any $\epsilon>0$,
there exists $C_{N, \epsilon, m,k} > 0$ such that
\begin{align}\label{f1}
 \begin{split}
& \;\theta\int_{\Omega} |u|^{\theta+1}\zeta^{4k} dx -\left(m-1\right)\left(1+ C_{m,k}\epsilon\right) \int_{\Omega}|\D u|^{m} \zeta^{4k} dx\\
\leq  & \;C_{N, \epsilon, m,k}\int_{\Omega} |u|^{m} \left(|\D \zeta|^{m}+|\nabla \zeta|^{2m} +|\nabla^{2} \zeta|^{m}\right)\zeta^{4k-2m} dx.
  \end{split}
\end{align}
Moreover, multiplying the equation \eqref{1} by $u \zeta^{4k}$ and
integrating by parts, there holds
\begin{align*}
 \int_{\Omega}|\D u|^{m} \zeta^{4k} dx-\int_{\Omega} |u|^{\theta+1}\zeta^{4k} dx \leq \int_{\Omega}|u||\Delta u|^{m-1}
  |\D(\zeta^{4k})| dx
  + C\int_{\Omega}|\Delta u|^{m-1} |\nabla u||\nabla (\zeta^{4k})| dx.
\end{align*}
Using Young's inequality and applying again Lemma \ref{l.2.7a}, we can conclude that for any $\epsilon>0$, there exists $C_{N,\epsilon, m,k} > 0$
 such that
\begin{align}\label{f12}
\begin{split}
& \;\left(1- C_{m,k}\epsilon\right) \int_{\Omega}|\D u|^{m} \zeta^{4k} dx-\int_{\Omega} |u|^{\theta+1}\zeta^{4k} dx\\
\leq& \; C_{N,\epsilon, m,k}\int_{\Omega} |u|^{m} \left(|\D \zeta|^{m}+|\nabla \zeta|^{2m} +|\nabla^{2} \zeta|^{m}\right)\zeta^{4k-2m} dx.
\end{split}
\end{align}

Taking $\epsilon > 0$ but small enough, multiplying \eqref{f12} by $\frac{(m-1)(1+ 2C_{m,k}\epsilon)}{1- C_{m,k}\epsilon}$, adding it with \eqref{f1}, we get
\begin{align*}
&(m-1) C_{m,k}\epsilon \int_{\Omega}|\D u|^{m} \zeta^{4k} dx + \left[\theta-\frac{(m-1)(1+ 2C_{m,k}\epsilon)}{1- C_{m, k}\epsilon}\right]\int_{\Omega} |u|^{\theta+1}\zeta^{4k} dx\\
\leq  & \; C_{N,\epsilon, m,k}\int_{\Omega} |u|^{m} \left(|\D \zeta|^{m}+|\nabla \zeta|^{2m} +|\nabla^{2} \zeta|^{m}\right)\zeta^{4k-2m} dx.
\end{align*}
As $\theta > m - 1 > 1$, using $\epsilon > 0$ small enough, we have
\begin{align}\label{f1xx2xx}
\int_{\Omega}|\D u|^{m} \zeta^{4k} dx +\int_{\Omega} |u|^{\theta+1}\zeta^{4k} dx\leq C\int_{\Omega} |u|^{m} \left(|\D \zeta|^{m}+|\nabla \zeta|^{2m} +|\nabla^{2} \zeta|^{m}\right)\zeta^{4k-2m} dx.
\end{align}
For $k\geq  \frac{m(\theta+1)}{2(\theta+1-m)}$ so that $4km\leq(4k-2m)(\theta+1),$ Applying H\"older inequality, we conclude then
\begin{align*}
& \;\int_{\Omega}|\D u|^{m} \zeta^{4k} dx +\int_{\Omega} |u|^{\theta+1}\zeta^{4k} dx\\
\leq  & \; C\left[\int_{\Omega}\left(|\D \zeta|^{m}+|\nabla \zeta|^{2m}
 +|\nabla^{2} \zeta|^{m}\right)^{\frac{\theta+1}{\theta-(m-1)}} dx\right]^{\frac{\theta-(m-1)}{\theta+1}}
 \left(\int_{\Omega} |u|^{\theta+1}\zeta^{\frac{(4k-2m)(\theta+1)}{m}}dx\right)^{\frac{m}{\theta+1}}\\
 \leq  & \; C\left[\int_{\Omega}\left(|\D \zeta|^{m}+|\nabla \zeta|^{2m}
 +|\nabla^{2} \zeta|^{m}\right)^{\frac{\theta+1}{\theta-(m-1)}} dx\right]^{\frac{\theta-(m-1)}{\theta+1}}
 \left(\int_{\Omega} |u|^{\theta+1}\zeta^{4k}dx\right)^{\frac{m}{\theta+1}}.
 \end{align*}
We get readily the estimate \eqref{new7}. \qed

\medskip
Now we choose $\phi_0$ a cut-off function in $C_c^\infty(B_2)$ verifying $0 \leq \phi_0 \leq 1$ and $\phi_0=1$ in $B_1$. Applying \eqref{new7} with $\zeta = \phi_0(R^{-1}x)$
for $R > 0$, there holds
 \begin{align*}
\int_{B_{R}} |u|^{\theta+1}dx\leq \int_{\R^N} |u|^{\theta+1}\zeta^{4k} dx \leq C R^{\frac{N(\theta -m+1)}{\theta + 1} -2m}.
\end{align*}
Under the assumption \eqref{new8}, tending $R \to \infty$, we obtain $u \equiv 0$, we prove then
Proposition \ref{p12bis}, hence the case $\theta p > 1 > p > 0$ for Proposition \ref{p12}.

\medskip
\noindent{\bf Proof of Lemma \ref{l.2.7a}.} A direct calculation gives
\begin{align}
  \label{new1}
  \begin{split}
 \int_{\Omega} |\nabla u|^m |\nabla \zeta|^{m}\zeta^{4k-m} dx = &\;\int_{\Omega} \nabla u\cdot\nabla u|\nabla u|^{m-2}|\nabla \zeta|^{m}\zeta^{4k-m}dx\\
 =& \; - \int_{\Omega}  \mathrm{div} \left(\nabla u|\nabla u|^{m-2}\right) u|\nabla \zeta|^{m}\zeta^{4k-m}dx\\
 &\; -\int_{\Omega} u|\nabla u|^{m-2}\nabla u \cdot\nabla\left(|\nabla \zeta|^{m}\zeta^{4k-m}\right)dx\\
 := & \; I_1 + I_2.
 \end{split}
\end{align}
The integral $I_1$ can be estimated as
 \begin{align*}
 I_1 & = -\left(m-2\right) \int_{\Omega}  u |\nabla u|^{m-4}|\nabla \zeta|^{m}\nabla^{2}u (\nabla u, \nabla u)\zeta^{4k-m}dx -\int_{\Omega}  u\D u|\nabla u|^{m-2}|\nabla \zeta|^{m}\zeta^{4k-m}dx &\;\\
 & \leq C_m\int_{\Omega}  |u|| \nabla^{2} u||\nabla u|^{m-2}|\nabla \zeta|^{m}\zeta^{4k-m}dx + \int_{\Omega}  |u||\D u||\nabla u|^{m-2}|\nabla \zeta|^{m}\zeta^{4k-m}dx.
\end{align*}
Applying Young's inequality, there holds, for any $\epsilon > 0$,
 \begin{align}\label{new3}
 \begin{split}
 & \int_{\Omega}  |u|| \D u||\nabla u|^{m-2}|\nabla \zeta|^{m}\zeta^{4k-m}dx\\
 \leq & \;C_{ \epsilon, m}\int_{\Omega}  |u|^{\frac{m}{2}}| \D u|^{\frac{m}{2}}|\nabla \zeta|^{m}\zeta^{4k-m}dx + \epsilon\int_{\Omega}|\nabla u|^{m}|\nabla \zeta|^{m}\zeta^{4k-m}dx\\
 \leq &\;C_{\epsilon, m}\int_{\Omega}  |u|^{m}|\nabla \zeta|^{2m}\zeta^{4k-2m}dx+ \epsilon\int_{\Omega}| \D u|^{m}\zeta^{4k}dx + \epsilon\int_{\Omega}|\nabla u|^{m}|\nabla \zeta|^{m}\zeta^{4k-m}dx.
\end{split}
\end{align}
On the other hand,
 \begin{align}\label{new2}
 \begin{split}
 & \;\int_{\Omega}  |u|| \nabla^{2} u||\nabla u|^{m-2}|\nabla \zeta|^{m-2+2}\zeta^{4k-m}dx\\
 \leq &\;  C_{\epsilon, m}\int_{\Omega}  |u|^{\frac{m}{2}}| \nabla^{2} u|^{\frac{m}{2}}|\nabla \zeta|^{m}\zeta^{4k-m}dx + \epsilon\int_{\Omega}|\nabla u|^{m}|\nabla \zeta|^{m}\zeta^{4k-m}dx\\
 \leq &\;C_{\epsilon, m}\int_{\Omega}  |u|^{m}|\nabla \zeta|^{2m}\zeta^{4k-2m}dx+ \epsilon\int_{\Omega}| \nabla^{2} u|^{m}\zeta^{4k}dx+ \epsilon\int_{\Omega}|\nabla u|^{m}|\nabla \zeta|^{m}\zeta^{4k-m}dx.
\end{split}
\end{align}
Now we shall estimate the integral
$$\int_{\Omega}| \nabla^{2} u|^{m}\zeta^{4k}dx.$$
Remark that there exists $C_0(N, m) > 0$ such that
\begin{align}
 \label{W2m}
 \int_{\R^N} |\nabla^2\varphi|^m dx \leq C_0(N, m)\int_{\R^N}|\Delta\varphi|^m dx, \quad \forall\; \varphi \in W^{2, m}(\R^N).
\end{align}
We can prove it firstly for $\varphi \in W^{2, m}_0(B_1)$ with elliptic theory, then for general $\varphi \in W^{2, m}(\R^N)$ with
approximation and scaling argument. As $u\zeta \in W^{2, m}_0(\Omega) \subset W^{2, m}(\R^N)$, \eqref{W2m} implies that
\begin{align*}
\int_{\Omega} |\nabla^2(u\zeta^{\frac{4k}{m}})|^m dx \leq & \;
C_0(N, m) \int_{\Omega}|\D (u \zeta^{\frac{4k}{m}})|^m dx\\
\leq & \; C_{N, m}\int_{\Omega}|\D u|^{m}\zeta^{4k} dx + C_{N, m,k}\int_{\Omega}|\nabla u|^{m}|\nabla \zeta|^{m}\zeta^{4k-m}dx\\
  &\; + C_{N, m,k}\int_{\Omega}|u|^{m}\left(|\nabla \zeta|^{2m} +|\nabla^{2} \zeta|^{m}\right)\zeta^{4k-2m}dx.
\end{align*}
Let $k > m$, we get then
 \begin{align}\label{new4}
 \begin{split}
 \int_{\Omega}| \nabla^{2} u|^{m}\zeta^{4k}dx
 \leq& \;C\int_{\Omega}|\nabla^2(u\zeta^{\frac{4k}{m}})|^m dx +C_{m,k}\int_{\Omega}|\nabla u|^{m}|\nabla \zeta|^{m}\zeta^{4k-m}dx\\
 &\;+ C_{m,k}\int_{\Omega}|u|^{m}\left(|\nabla \zeta|^{2m} +|\nabla^{2} \zeta|^{m}\right)\zeta^{4k-2m}dx\\
 \leq & \; C_{N,m}\int_{\Omega}|\D u|^{m}\zeta^{4k} dx + C_{N, m,k}\int_{\Omega}|\nabla u|^{m}|\nabla \zeta|^{m}\zeta^{4k-m}dx\\
  &\; + C_{N, m,k}\int_{\Omega}|u|^{m}\left(|\nabla \zeta|^{2m} +|\nabla^{2} \zeta|^{m}\right)\zeta^{4k-2m}dx.
 \end{split}
\end{align}
Combining \eqref{new3}, \eqref{new2} and \eqref{new4}, we arrive at
 \begin{align}\label{SSS0.255}
 \begin{split}
I_1 \leq &\;C_{N,m,k}\epsilon \int_{\Omega}|\D u|^{m}\zeta^{4k} dx +C_m\epsilon\int_{\Omega}|\nabla u|^{m}|\nabla \zeta|^{m}\zeta^{4k-m}dx \\
  &+\; C_{N,\epsilon, m, k} \int_{\Omega}|u|^{m}\left(|\nabla \zeta|^{2m} +|\nabla^{2} \zeta|^{m}\right)\zeta^{4k-2m}dx.
 \end{split}
\end{align}

Furthermore, by Young’s inequality,
 \begin{align}\label{SS0.255}
 \begin{split}
I_2 = & - m\int_{\Omega} u|\nabla u|^{m-2} |\nabla \zeta|^{m-2}\nabla^2\zeta(\nabla\zeta, \nabla u)\zeta^{4k-m} dx\\
& -(4k-m)\int_{\Omega} u|\nabla u|^{m-2}|\nabla \zeta|^{m}(\nabla u\cdot \nabla\zeta)\zeta^{4k-m-1}dx\\
\leq &\; C_{ m,k} \int_{\Omega}|u||\nabla u|^{m-1}|\nabla \zeta|^{m-1}\left(|\nabla \zeta|^{2} +|\nabla^{2} \zeta|\right)\zeta^{4k-m-1}dx\\
\leq &\; C_{\epsilon, m,k} \int_{\Omega}|u|^{m}\left(|\nabla \zeta|^{2m} +|\nabla^{2} \zeta|^{m}\right)\zeta^{4k-2m}dx
+ \epsilon\int_{\Omega}|\nabla u|^{m}|\nabla \zeta|^{m}\zeta^{4k-m}dx.
 \end{split}
\end{align}
Combining \eqref{SSS0.255}--\eqref{SS0.255} with \eqref{new1}, one concludes
 \begin{align*}
(1-C_{N, m,k}\epsilon)\int_{\Omega} |\nabla u|^m |\nabla \zeta|^{m}\zeta^{4k-m} dx
& \leq C_{N,\epsilon, m ,k} \int_{\Omega}|u|^{m}\left(|\nabla \zeta|^{2m} +|\nabla^{2} \zeta|^{m}\right)\zeta^{4k-2m}dx\\
  &+C_{N, m, k}\epsilon\int_{\Omega}|\D u|^{m}\zeta^{4k} dx.
\end{align*}
This means that \eqref{new13} holds true for $\epsilon > 0$ small enough, hence for any $\epsilon > 0$.\qed

\section{Proof of Theorem \ref{main4}.}
\setcounter{equation}{0}
In this section, we prove Theorem \ref{main4}. As already mentioned, we need only to handle the case $p\theta > 1$. We use first the classification for stable solutions, Proposition \ref{p12} to obtain the decay estimates for stable at infinity solutions of \eqref{1.1}.

\begin{lem}\label{newla}
 Let $p, \theta > 0$ verify $p\theta > 1$ and \eqref{LEbis}. Let $(u,v)$ be a solution of \eqref{1.1}  which is stable outside a compact set.
 Then there exists a constant C such that
\begin{align}\label{new5}
\sum_{k\leq 2}\Big[ |x|^{\alpha+k} |\nabla^k u(x)|+|x|^{\beta+k} |\nabla^k v(x)|\Big]\leq C,\quad  \forall\; x \in \R^N.
\end{align}
\end{lem}
\textbf{Proof.} Assume that $(u, v)$ is stable outside $B_{R_0}$. Denote $$W(x)=\sum_{k\leq 2}\left[ |\nabla^k u(x)|^{{\frac{1}{\alpha+k}}}+ |\nabla^k v(x)|^{\frac{1}{\beta+k}}\right].$$
Suppose that \eqref{new5} does not hold true. Let $d(x) = \|x\|-R_0$, there holds
$$\sup_{\R^N\backslash B_{R_0}}W(x)d(x) =\infty,$$
or equally there exists a sequence $(x_n)$ such that $\|x_n\| > R_0$ and $W(x_{n})d(x_{n})>n$ for $n \geq 1$. Since $(u, v)$ are smooth in $\mathbb{R}^N$, then $d(x_{n})\rightarrow \infty$. By the doubling lemma \cite{pqs}, there exists another sequence $(y_{n})$  such that for any $n \geq 1$, $\|y_n\| > R_0$,
\begin{enumerate}
  \item [$(i)$] $W(y_{n})d(y_{n})\geq n$;
  \item [$(ii)$] $W(y_{n})\geq W(x_{n})$;
  \item [$(iii)$] $W(z)\leq 2 W(y_{n})$  for  $|z|> R_0$  such that  $|z-y_{n}|\leq \frac{n}{W(y_{n})}.$
\end{enumerate}

 Let $(u,v)$ be a solution of \eqref{1.1}, consider the sequence of functions
 \begin{align*}
\wt u_{n}(x) = \lambda_{n}^{\alpha}u(y_n + \lambda_{n}x), \quad \wt v_{n}(x) = \lambda_{n}^{\beta} v(y_n + \lambda_{n}x), \quad \mbox{with }\; \lambda_{n}=W(y_{n})^{-1}.
\end{align*}
It's well known that $(\wt u_n, \wt v_n)$ are a sequence of solutions to \eqref{1.1}. Moreover,
$$W_n(x) := \sum_{k\leq 2}\left( |\nabla^k \wt u_n(x)|^{\frac{1}{\alpha+k}}+ |\nabla^k \wt v_n(x)|^{{\frac{1}{\beta+k}}}\right) = \lambda_n W(y_n + \lambda_n x),
\quad \forall\; x \in \R^N.$$
By $(i)$, we have $B_{n\lambda_n}(y_{n}) \subset \R^N\backslash B_{R_0}$, and we can readily check that $(\wt u_n, \wt v_n)$ is stable in $B_{n}$ since $(u,v)$ is stable in $\R^N\backslash B_{R_0}$. Using $(iii)$, there holds, for all $n \geq 1$,
\begin{align}\label{newAnewest5}
W_n(x) \leq 2W_n(0) = 2\quad \mbox{in }\; B_n.
\end{align}
From \eqref{newAnewest5} and standard elliptic theory, up to a subsequence, $(\wt u_n, \wt v_n)$ converges to $(u_\infty, v_\infty)$ in
 $C^2_{loc}({\mathbb R}^N)$. Therefore
$$\sum_{k\leq 2}\left( |\nabla^k u_\infty(0)|^{\frac{1}{\alpha+k}}+ |\nabla^k v_\infty(0)|^{{\frac{1}{\beta+k}}}\right)  = 1.$$
So $(u_{\infty}, v_\infty)$  is nontrivial.  Clearly, $(u_\infty, v_\infty)$ a smooth positive solution to \eqref{1.1}. Using again the elliptic theory, it's not difficult to see that $(u_\infty, v_\infty)$ is stable in $\R^N$. However, this contradicts Proposition \ref{p12} since $p, \theta$ verifies \eqref{LE}. Hence the hypothesis was wrong, i.e.~the estimate \eqref{new5} holds true.  \qed

\medskip
Another tool is the following classical Pohozaev identity (see \cite{em, PJ, sou}).
\begin{lem}\label{newlem1.1}
Let $(u,v)$ be a solution to \eqref{1.1}. Therefore for any regular bounded domain $\Omega$,
\begin{align}\label{new6}
\begin{split}
&\; \frac{2(p+1)- pN}{p+1}\int_{\Omega} v^{p+1} dx +\frac{N}{\theta+1} \int_{\O} u^{\theta+1} dx \\
=& \int_{\partial\Omega}u^{\theta+1}(\nu\cdot x)d\sigma-\frac{p}{p+1} \int_{\partial\Omega}v^{p+1}(\nu\cdot x)d\sigma
+\int_{\partial\Omega}\frac{\p v}{\p \nu}(\nabla u\cdot x)d\sigma-\int_{\partial\Omega}v\frac{\p(\nabla u\cdot x)}{\p \nu}d\sigma.
\end{split}
 \end{align}
 \end{lem}

We claim then
 \begin{lem}\label{newl.552.GR7a}
Let $p, \theta > 0$ satisfy $p\theta > 1$ and \eqref{LEbis}. If $(u,v)$ be a solution of \eqref{1.1} which is stable outside a compact set, then $v \in L^{p+1}(\mathbb{R}^N)$, $u\in L^{\theta+1}(\mathbb{R}^N)$ and
\begin{align}\label{new8newest0}
\frac{2(p+1)-pN }{p+1}\int_{\mathbb{R}^N}v^{p+1}dx + \frac{N}{\theta+1}\int_{\mathbb{R}^N}u^{\theta+1}dx = 0.
 \end{align}
\end{lem}

\noindent{\bf Proof .}  By \eqref{new5}, we  have (noticing that $\alpha(\theta+1)=(p+1)\beta=2+\alpha+\beta$)
\begin{align*}
u^{\theta+1}(x)+v^{p+1}(x) \leq C\left(1+ |x|\right)^{-(2 + \alpha+\beta)} \quad \mbox{in }\; \R^N.
 \end{align*}
By \eqref{LEbis}, then $v \in L^{p+1}(\mathbb{R}^N)$, $u\in L^{\theta+1}(\mathbb{R}^N) .$ Using  Lemma \ref{newlem1.1} with $\Omega =B_{R}$, we deduce that
\begin{align}\label{newac2k.5}
\begin{split}
& \frac{2(p+1)-pN }{p+1}\int_{B_{R}}v^{p+1}dx+ \frac{N}{\theta+1}\int_{B_{R}}|u|^{\theta+1}dx\\
=&\;\int_{\partial B_{R}}\left[R\frac{\partial  u}{\partial r}\frac{\partial v}{\partial r}-\frac{pR}{p+1}v^{p+1}
 -v\frac{\partial (\nabla u\cdot x)}{\partial r}+\frac{R}{\theta+1}|u|^{\theta+1}\right]d\sigma.
\end{split}
 \end{align}
Using again \eqref{new5} and $N < 2 + \alpha+\beta$, we deduce that
 $$ \int_{\partial B_{R}}\left[R\frac{\partial  u}{\partial r}\frac{\partial v}{\partial r}-\frac{R}{p+1}v^{p+1}
 -v\frac{\partial (\nabla u\cdot x)}{\partial r}+\frac{R}{\theta+1}|u|^{\theta+1}\right]d\sigma\rightarrow 0, \quad \mbox{as }\; R\rightarrow \infty.$$
Taking the limit $R \to\infty$ in \eqref{newac2k.5}, the claim follows. \qed

\medskip\noindent
{\bf Proof of Theorem \ref{main4} completed.} We are now in position to conclude. Suppose that $(u, v)$ is a solution to \eqref{1.1} stable at infinity with $p\theta > 1$ verifying \eqref{LEbis}. Choose $\phi_0$ a cut-off function in $C_c^\infty(B_2)$ verifying $0 \leq \phi_0 \leq 1$ and $\phi_0=1$ in $B_1$. Denote $\zeta = \phi_0(R^{-1}x)$ and $A_R= B_{2R}\backslash B_R$. By the system \eqref{1.1}, there holds
\begin{align*}
 \int_{B_{2R}}v^{p+1} \zeta dx-\int_{B_{2R}} u^{\theta+1} \zeta dx & =  \int_{B_{2R}}u\zeta \Delta v dx - \int_{B_{2R}} u\zeta \Delta v dx\\
& = \int_{B_{2R}}v\Big(2\nabla u\cdot \nabla \zeta + u\Delta\zeta\Big) dx\\
& \leq \frac{C}{R^{2}} \int_{A_R}uv dx
  + \frac{C}{R} \int_{A_R}v |\nabla u| dx.
\end{align*}
Using \eqref{new5}, and tending $R\rightarrow \infty$, as $N < 2 + \alpha+\beta,$ we have
$$\int_{\mathbb{R}^N}v^{p+1} dx=\int_{\mathbb{R}^N}u^{\theta+1} dx.$$
Substituting this in \eqref{new8newest0},
\begin{align*}
\left(\frac{2(p+1) - pN}{p+1}+\frac{N}{\theta+1}\right)\int_{\mathbb{R}^N}u^{\theta+1}dx=0.
\end{align*}
As \eqref{LEbis} implies that
$$\frac{2(p+1) - pN}{p+1}+\frac{N}{\theta+1} = 2 - \frac{(p\theta - 1)N}{(p+1)(\theta + 1)} = 2 - \frac{2N}{2+\alpha + \beta} > 0,$$
$u \equiv 0$ in $\R^N$ which is absurd, so we are done.\qed

\section{\bf Proof of Theorem \ref{main3}.}
\setcounter{equation}{0}
The approach is similar to that for Theorem \ref{main4}. We derive first some integral estimates thanks to Lemma \ref{lemnewBN}. Suppose that $u$ is stable outside the ball $B_{R_0}$. Let $R > R_0 + 3$ and $\zeta \in C^2_c(\R^N\backslash B_{R_0})$ verifying that $0 \leq \zeta \leq 1$ and
$$\zeta(x) =\left\{
\begin{array}{ll}
       0 \quad \mbox{for}\; \|x\|\leq R_{0}+1,\; \|x\|\geq 2R,\\
      1\quad \mbox{for}\; R_{0}+2 \leq \|x\| \leq R.
       \end{array}
     \right.$$
Clearly, we can assume that there exists $C > 0$ independent on $R$ such that
$$\|\zeta\|_{C^2(B_{R_0+2})} \leq C \quad \mbox{and} \quad R|\nabla \zeta(x)| + R^2|\nabla^2 \zeta(x)|\leq C \;\; \mbox{in } A_R = B_{2R}\backslash B_R.$$
Applying the estimate \eqref{new7} with $\zeta$, we get readily
\begin{align}\label{AvxZK}
 \begin{split}
& \; \int_{R_{0}+2\leq \|x\|\leq R} |\D u|^m dx + \int_{R_{0}+2\leq \|x\|\leq R} |u|^{\theta+1} dx \leq C\left(1 + R^{N-\frac{2m(\theta + 1)}{\theta -(m-1)}}\right).
 \end{split}
\end{align}
Using \eqref{new8} and tending $R \to \infty$, we have then
\begin{align}
\label{new11}
u \in L^{\theta + 1}(\R^N) \quad \mbox{and}\quad \Delta u \in L^m(\R^N).
\end{align}
By H\"older's inequality, there holds
\begin{align*}
R^{-2m} \int_{B_{R}} |u|^mdx  \leq CR^{\frac{N(\theta+1-m)}{\theta+1}-2m}\left(\int_{B_R} |u|^{\theta+1}dx\right)^{\frac{m}{\theta+1}}.
\end{align*}
On the other hand, by standard scaling argument, there exists $C > 0$ such that for any $R > 0$, any $u \in W^{2, m}(A_R)$ with $A_R = B_{2R}\backslash B_R$,
\begin{align*}
R^{-m} \int_{A_{R}} |\nabla u|^mdx  \leq C \int_{A_{R}} |\Delta u|^mdx  + CR^{-2m} \int_{A_{R}} |u|^mdx.
\end{align*}
Therefore, under the assumptions of Theorem \ref{main3}, we get
\begin{align}
\label{new9}
R^{-2m} \int_{A_{R}} |u|^mdx + R^{-m} \int_{A_{R}} |\nabla u|^mdx \to 0 \quad \mbox{as }\; R \to \infty.
\end{align}
Let $\zeta(x) = \phi_0(R^{-1}x)$ with a standard cut-off function $\phi_0 \in C_c^2(B_2)$, $\phi_0\equiv 1$ in $B_1$. Applying the estimate \eqref{new4} and using \eqref{new11}--\eqref{new9}, there holds
\begin{align}
\label{new10}
\int_{R^N}|\nabla^2 u|^m dx < \infty.
\end{align}

However, as we have mentioned, the weak solutions of \eqref{1} are in general not belongs to $C^2$, so we cannot use the standard Pohozaev identity similar to \eqref{new6} because of the boundary terms. We show here a variant of the Pohozaev identity, which proof is given in the appendix for the convenience of the readers.
\begin{lem}\label{lem1h.1}
Let $u$ be a  weak solution to \eqref{1} with $m > 2$. Then for any $\psi \in C_c^{2}(\Omega)$,
\begin{align}\label{new12}
\begin{split}
&\; \frac{N}{\theta+1}\int_{\Omega} |u|^{\theta+1}\psi dx
- \frac{N-2m}{m}\int_{\Omega} |\Delta u|^{m}\psi dx \\
 =& \;-\frac{1}{\theta+1}\int_{\Omega} |u|^{\theta+1}(\nabla \psi\cdot x) dx + \frac{1}{m}\int_{\Omega} (\nabla\psi\cdot x)|\Delta u|^{m} dx\\
& \; - \int_{\Omega} |\Delta u|^{m-2}\Big[2\Delta u(\nabla u\cdot \nabla \psi)  + 2\Delta u \nabla^2(x\cdot\nabla\psi) + \Delta u(\nabla u\cdot x)\Delta \psi \Big]dx.
\end{split}
 \end{align}
 \end{lem}

This implies that if u is a  weak solution of \eqref{1}, stable at infinity with $1 < m-1<\theta$ and $N$ verifying \eqref{new8}, then
\begin{align}\label{8newest0}
 \frac{N-2m }{m}\int_{\mathbb{R}^N}|\Delta u|^{m}dx=\frac{N}{\theta+1}\int_{\mathbb{R}^N}|u|^{\theta+1}dx.
 \end{align}
Indeed, let $\psi$ in \eqref{new12} be defined by $\psi(x) = \phi_0(R^{-1}x)$ with a standard cut-off function $\phi_0 \in C_c^2(B_2)$, $\phi_0\equiv 1$ in $B_1$.
Denote the right hand side in \eqref{new12} by $I_R$. Remark that $\nabla \psi \ne 0$ only in $A_R = B_{2R}\backslash B_R$ and $\|\nabla^k \psi\|_\infty \leq C_k R^{-k}$, we obtain readily
\begin{align*}
|I_R| \leq C\int_{A_R} \Big(|u|^{\theta+1} + |\Delta u|^m\Big) dx + \frac{C}{R}\int_{A_R}|\Delta u|^{m-1}|\nabla u|dx + C\int_{A_R}|\Delta u|^{m-1}|\nabla^2 u| dx
\end{align*}
Thanks to the estimates \eqref{new11}-\eqref{new10} and H\"older's inequality, clearly $\lim_{R\to \infty}I_R = 0$, hence we get \eqref{8newest0}.

\medskip
On the other hand, using $u \psi$ as test function in \eqref{1} , we get
\begin{align*}
 \int_{B_{2R}}|\Delta u|^{m} \psi dx-\int_{B_{2R}} |u|^{\theta+1} \psi dx &\leq C\int_{B_{2R}}|u||\Delta u|^{m-1}
  |\D \psi| dx
  + C\int_{B_{2R}}|\Delta u|^{m-1} |\nabla u||\nabla\psi| dx\\
& \leq \frac{C}{R^2}\int_{A_R}|u||\Delta u|^{m-1}dx
  + \frac{C}{R}\int_{A_R}|\Delta u|^{m-1}|\nabla u|dx.
\end{align*}
Apply H\"older's inequality, \eqref{new11}--\eqref{new9} and tending $R$ to $\infty$,  we obtain
 \begin{align}\label{llfou}
\int_{\mathbb{R}^N}|u|^{\theta+1}dx=\int_{\mathbb{R}^N} |\D u|^{m}dx.
\end{align}
Combining \eqref{8newest0} and \eqref{llfou}, one obtains
\begin{align*}
\left(\frac{N-2m }{m}-\frac{N}{\theta+1}\right)\int_{\mathbb{R}^N}|u|^{\theta+1}dx=0.
\end{align*}
We are done, since \eqref{new8} implies that $\frac{N-2m }{m}-\frac{N}{\theta+1} < 0$.\qed

\section{Appendix }
\setcounter{equation}{0}
We prove here the Lemma \ref{lem1h.1}. Let $\psi\in C_c^{2}(\Omega)$, multiplying equation \eqref{1} by $\nabla u\cdot x \psi$ and integrating by parts, we get
\begin{align*}
& \int_{\Omega}|u|^{\theta-1}u (\nabla u\cdot x)\psi dx\\ = & \;  \int_{\Omega}|\Delta u|^{m-2}\Delta u\D(\nabla u\cdot x\psi)dx\\
= & \; \int_{\Omega}\; |\Delta u|^{m-2}\Delta u \Big[(\nabla(\Delta u)\cdot x)\psi+2\Delta u \psi+2\nabla(\nabla u\cdot x)\cdot\nabla \psi+(\nabla u\cdot x)\Delta\psi \Big]dx.
\end{align*}
Direct calculation yields $\nabla(\nabla u\cdot x)\cdot\nabla \psi = \nabla^{2}u (x, \nabla \psi) + (\nabla u\cdot \nabla \psi)$ and
\begin{align*}
 \int_{\Omega}\; |\Delta u|^{m-2}\Delta u \Big[(\nabla(\Delta u)\cdot x)\psi+2\Delta u \psi\Big]dx
   &= \int_{\Omega} \frac{\nabla |\Delta u|^m}{m}\cdot x \psi dx+ 2\int_{\Omega}|\Delta u|^{m}\psi dx\\
   &=\frac{2m-N }{m}\int_{ \Omega}|\Delta u|^m\psi dx-\frac{1}{m}\int_{\Omega}|\Delta u|^m(\nabla \psi\cdot x)dx.
\end{align*}
Moreover,
\begin{align*}
\int_\Omega |u|^{\theta-1}u (\nabla u\cdot x)\psi
  dx & = - \frac{1}{\theta + 1}\int_\Omega |u|^{\theta +1} {\rm div}(\psi x)dx \\
&=-\frac{N}{\theta+1} \int_{\Omega}|u|^{\theta+1} \psi dx - \frac{1}{\theta+1} \int_{\Omega}|u|^{\theta+1} x\cdot \nabla \psi dx.
\end{align*}
Therefore, the claim follows by regrouping the above equalities. \qed
\vskip 1cm

\section*{Acknowledgment}
The authors are partly supported by the CNRS-DGRST Project No.~EDC26348.

\end{document}